\documentclass[10pt]{article}
\usepackage{amsthm,amsmath,amssymb,xcolor,stmaryrd}
\usepackage{authblk}
\usepackage{graphicx}
\usepackage{enumitem}
\usepackage{amsmath}
\usepackage{amssymb}
\usepackage{multicol}
\usepackage{xcolor}

\newcommand{\ignore}[1]{}

\newcommand{\epfun}[1]{\check \ve\left(#1\right)}
\newcommand{\IG}{{\rm G}^{\text{\sc a}}}
\newcommand{\BG}{{\rm G}^{\text{\sc b}}}
\newcommand{\CG}{{\rm G}^{\text{\sc c}}}
\newcommand{\io}{o^{\text{\sc a}}}
\newcommand{\bo}{o^{\text{\sc b}}}
\newcommand{\co}{o^{\text{\sc c}}}
\newcommand{\fsi}[2]{#1\llbracket #2 \rrbracket }

\newcommand{\bchk}{( k^+)}

\newcommand{\bchkb}{\langle k^+\rangle}

\newcommand{\bchkc}{\{ k^+\}}

\newcommand{\bch}[2]{ 
\left ( \begin{smallmatrix}
  #2\\
  #1
\end{smallmatrix}
\right )
}

\newcommand{\bchb}[2]{ 
 \langle \begin{smallmatrix}
  #2\\
  #1
\end{smallmatrix}
\rangle
}

\newcommand{\bcc}[2]{ 
\left \{ \begin{smallmatrix}
  #2\\
  #1
\end{smallmatrix}
\right \}
}

\newcommand{\gnf}{=_k}
\newcommand{\gnfp}{=_{k+1}}

\newcommand{\gnfm}{\equiv_k}

\newcommand{\pk}[1]{#1\bchb k\omega}
\newcommand{\zk}[1]{#1\bch k\omega}

\newcommand{\pskpe}[1]{#1\bchb {k+1}\omega}
\newcommand{\zskpe}[1]{#1\bch {k+1}\omega}

\newcommand{\ck}[1]{#1 \bcc k\omega}
\newcommand{\cskpe}[1]{#1\bcc{k+1}\omega}
\newcommand{\ve}{\varepsilon}

\renewcommand{\phi}{{\overline{\varphi}}}

\newcommand{\al}{{\alpha}}

\newcommand{\be}{{\beta}}

\newcommand{\ga}{\gamma}

\newcommand{\la}{\lambda}

\newcommand{\om}{\omega}

\newtheorem{theorem}{Theorem}[section]
\newtheorem{lemma}[theorem]{Lemma}
\newtheorem{definition}[theorem]{Definition}

\newtheorem{proposition}[theorem]{Proposition}

\begin{document}
\title{Intermediate Goodstein principles}
\author{David Fern\'andez-Duque\thanks{David.FernandezDuque@UGent.be} }
\author{Oriola Gjetaj\thanks{Oriola.Gjetaj@UGent.be}}
\author{Andreas Weiermann\thanks{Andreas.Weiermann@UGent.be}}
\affil{Department of Mathematics: Analysis, Logic and Discrete Mathematics, Ghent University}
\maketitle              
\begin{abstract}
The original Goodstein process proceeds by writing natural numbers in nested exponential $k$-normal form, then successively raising the base to $k+1$ and subtracting one from the end result.
Such sequences always reach zero, but this fact is unprovable in Peano arithmetic.
In this paper we instead consider notations for natural numbers based on the Ackermann function.
We define three new Goodstein processes, obtaining new independence results for $ {\sf ACA}_0$, ${\sf ACA}_0'$ and ${\sf ACA}_0^+$, theories of second order arithmetic related to the existence of Turing jumps.
\end{abstract}

\section{Introduction}

Goodstei{n'}s principle \cite{Goodsteinb} is arguably the oldest example of a purely number-theoretic statement known to be independent of $\sf PA$, as it does not require the coding of metamathematical notions such as G\"odel's provability predicate \cite{Godel1931}.
The proof proceeds by transfinite induction up to the ordinal $\varepsilon_0$ \cite{Goodstein1944}.
$\sf PA$ does not prove such transfinite induction, and indeed Kirby and Paris later showed that Goodstei{n'}s principle is unprovable in $\sf PA$ \cite{Kirby}.

Goodstei{n'}s original principle involves the termination of certain sequences of numbers.
Say that $m$ is in {\em nested (exponential) base-$k$ normal form} if it is written in standard exponential base $k$, with each exponent written in turn in base $k$.
Thus for example, $20 $ would become $2^{2^2}+2^2$ in nested base-$2$ normal form.
Then, define a sequence $(g_k(0))_{m\in\mathbb N}$ by setting $g_0(m) = m$ and defining $g_{k+1}(m)$ recursively by writing $ g_k(m)$ in nested base-$(k+2)$ normal form, replacing every occurrence of $k+2$ by $k+3$, then subtracting one (unless $g_k(m)=0$, in which case $g_{k+1}(m) =0$).

In the case that $m=20$, we obtain
\begin{align*}
g_0(20) & =20 = 2^{2^2}+2^2\\
g_1(20) &  = 3^{3^3}+3^3 -1 = 3^{3^3}+3^2\cdot 2 + 3\cdot 2 + 2\\
g_2(20) & = 4^{4^4}+4^2\cdot 2 + 4\cdot 2 +2 -1 = 4^{4^4}+4^2\cdot 2 + 4\cdot 2 +1,
\end{align*} 
and so forth.
At first glance, these numbers seem to grow superexponentially.
It should thus be a surprise that, as Goodstein showed, for every $m$ there is $k^*$ for which $g_{k^*}(m)=0$. 

By coding finite Goodstein sequences as natural numbers in a standard way, Goodstei{n'}s principle can be formalized in the language of arithmetic, but this formalized statement is unprovable in $\sf PA$. Independence can be shown by proving that the Goodstein process takes at least as long as stepping down the {\em fundamental sequences} below $\varepsilon_0$; these are canonical sequences $( \xi [n])_{n<\omega}$ such that $ \xi [n]< \xi$ for all $\xi$ and for limit $\xi$, $ \xi [n] \to \xi$ as $n\to \infty$.
For standard fundamental sequences below $\varepsilon_0$, $\sf PA$ does not prove that the sequence $\xi >  \xi [1] > \xi[1][2] > \xi[1][2][3] \ldots $ is finite.\smallskip

Exponential notation is not suitable for writing very big numbers (e.g.~Graham's number \cite{graham1969}), in which case it may be convenient to use systems of notation which employ faster-growing functions.
In \cite{FSGoodstein}, T.~Arai, S.~Wainer and the authors have shown that the Ackermann function may be used to write natural numbers, giving rise to a new Goodstein process which is independent of the theory ${\sf ATR}_0$ of {\em arithmetical transfinite recursion;} this is a theory in the language of second order arithmetic which is much more powerful than $\sf PA$.
The main axiom of ${\sf ATR}_0$ states that for any set $X$ and ordinal $\alpha$, the $\alpha$-Turing jump of $X$ exists; we refer the reader to \cite{Simpson:2009:SubsystemsOfSecondOrderArithmetic} for details.

The idea is, for each $k\geq 2$, to define a notion of Ackermannian normal form for each $m\in \mathbb N$.
Having done this, we can define Ackermannian Goodstein sequences analogously to Goodstei{n'}s original version.
The normal forms used in \cite{FSGoodstein} are defined using an elaborate `sandwiching' procedure first introduced in \cite{Singapore}, approximating a number $m$ by successive branches of the Ackermann function.
In this paper, we consider simpler, and arguably more intuitive, normal forms, also based on the Ackermann function.
By making variations of how these normal forms are treated, we show that these give rise to three different Goodstein-like processes, independent of ${\sf ACA}_0$, ${\sf ACA}'_0$ and ${\sf ACA}^+_0$, respectively.
As was the case for ${\sf ATR}_0$, these are theories of second order arithmetic which state that certain Turing jumps exist.
Recall that ${\sf ACA}_0$ is a theory of second-order arithmetic whose characteristic axiom is that if $X$ is a set of natural numbers, then the Turing jump of $X$ also exists as a set.
The more powerful theory ${\sf ACA}'_0$ asserts that, for all $n\in\mathbb N$ and $X\subseteq \mathbb N$, the $n$-Turing jump of $X$ exists, while ${\sf ACA}^+_0$ asserts that its $\omega$-jump exists; see \cite{Simpson:2009:SubsystemsOfSecondOrderArithmetic} for details.
The theory ${\sf ACA}_0$ is a conservative extension of Peano arithmetic, hence has proof-theoretic ordinal $\varepsilon_0$ \cite{Pohlers:2009:PTBook}.
The proof-theoretic ordinal of ${\sf ACA}'_0$ is $\varepsilon_\omega$ \cite{Afshari12}, and that of ${\sf ACA}^+_0$ is $\varphi_2(0)$ \cite{Leigh10}; we will briefly review these ordinals later in the text, but refer the reader to standard texts such as \cite{Pohlers:2009:PTBook,Schutte1977} for a more detailed treatment of proof-theoretic ordinals.
Preliminary versions of some of the results reported here appeared originally in \cite{GoodsteinCiE}.

\section{Basic definitions}

Our Goodstein processes will be based on a parametrized version of the Ackermann function, as given by the following definition.

\begin{definition}
For $a,b\in\mathbb N$ and $k\geq 2$ we define $A_a(k,b)$ by the following recursion.
As a shorthand, we write $A_ab$ instead of $A_a(k,b) $.
Then, we define
\begin{enumerate}
\item $A_0 b :=k^b$,
\item $A_{a+1} 0  :=A_a^k  0$,
\item $A_{a+1}(b+1):=A_a^k A_{a+1}b$.
\end{enumerate}
\end{definition}

Throughout the text we maintain the convention of writing $A_ab$ instead of $A_a(k,b) $ when $k\geq 2 $ is clear from context.
It is well known that for every fixed $a$, the function $b\mapsto A_a b $ is primitive recursive and the function
$a\mapsto A_a 0 $ is not primitive recursive.
It will be convenient to establish some simple lower bounds for the Ackermann function.

\begin{lemma}\label{lemmUseful}
Fix $k\geq 2$ and write $A_xy$ for $A_x(k,y)$.
Let $a,b\in\mathbb N$. Then,
\begin{enumerate}

\item If $a<{a'}$ then $A_ab\leq A_{{a'}}b$ and if $b<{b'}$ then $A_ab\leq A_{a}{b'}$.

\item\label{itMult} $k A_ab \leq A_a(b+1)\leq A_{a+1}b$.

\item\label{itSquare} If $k\geq 3$,
$A_a b \geq  b^2+b$.

\end{enumerate}
\end{lemma}

\begin{proof}
Most inequalities follow by induciton on $a$ with a secondary induction on $b$.
For \ref{itSquare}, we note that
$A_a b\geq A_0b\geq 3^b\geq b^2+b.$
\end{proof}

We use the Ackermann function to define $k$ normal forms for natural numbers.
These normal forms emerged from discussions with Toshiyasu Arai and Stan Wainer, 
which finally led to the definition of a more powerful normal form defined in \cite{Singapore} and used to prove termination in \cite{FSGoodstein}.

\begin{lemma}
Fix $k\geq 2$.
For all $m>0$, there exist unique $a,b, c \in \mathbb N$ such that
\begin{enumerate}
\item $m=A_a b+c$,
\item\label{itaMax} $A_a 0 \leq m <A_{a+1} 0 $,
\item $A_ab \leq m <A_{a}( b+1)$
\end{enumerate}

\end{lemma}

We write $m\gnf A_ab + c$ in this case.
This means that we have in mind an underlying context fixed by $k$ and that for the number $m$
we have uniquely associated the numbers $a,b,c$. Note that it could be possible that
$A_{a+1} 0 =A_a b $ for some $a,b$, so that we have to choose the right representation for the context; in this case, item \ref{itaMax} guarantees that $a$ is chosen to take the maximal possible value.

By rewriting all of $a,b,c$ in such a normal form, we obtain the {\em nested Ackermann $k$-normal form} of $m$.
If we only rewrite $a$ and $c$ in such a normal form (but not $b$), we arrive at the {\em index-nested Ackermann $k$-normal form} of $m$; note that in this case, $b$ should be regarded as a constant.
If we instead rewrite $b$ and $c$ iteratively, we arrive at the {\em argument-nested Ackermann $k$-normal form} of $ m $.
The following properties of normal forms are not hard to prove from the definitions.

\begin{lemma}\label{NF}
Fix $k\geq 2$.
\begin{enumerate}
\item If  $m =A_a b  + c + 1 $ is in $k$-normal form then  $ A_ab+c$ is in $k$-normal form as well.

\item $A_a^\ell 0 $ is in $k$-normal form for every $\ell$ such that $0<\ell<k$.

\item if $A_a b $ is in $k$-normal form, then for every $\ell<b$, the number $A_a \ell $ is also in $k$-normal form.
\end{enumerate}
\end{lemma}

Note that if $m\gnf A_ab+c$, it may still be the case that $c>A_ab$.
In this case, it will sometimes be convenient to count the number of occurrences of $A_ab$.
We thus write $m\gnfm A_ab\cdot q+c'$ if $m \gnf A_ab + c$ for $c= A_ab\cdot(q-1)+c'$ and $c'<A_ab$, and say that $A_ab\cdot q+c'$ is in {\bf extended normal form.}

\section{Proof-theoretic ordinals}

We work with standard notations for ordinals. We use the function $\xi\mapsto\varepsilon_\xi$ to enumerate the fixed points
of $\xi\mapsto \om^\xi$. With $\al,\be\mapsto \varphi_\al(\be)$ we denote the binary Veblen function, where $\beta\mapsto\varphi_\alpha(\beta)$ enumerates the common fixed points of all $\varphi_{\alpha'}$ with $\alpha'<\alpha$.
We often omit parentheses and simply write $\varphi_\alpha\beta$.
Then $\varphi_0\xi=\om^\xi$, $\varphi_1\xi=\varepsilon_\xi$, $\varphi_20$ is the first fixed point of the function $\xi\mapsto \varphi_1\xi$,
$\varphi_\omega 0$ is the first common fixed point of the function $\xi\mapsto \varphi_n \xi$, and $\Gamma_0$ is the first ordinal closed under
$\al,\be\mapsto \varphi_\al\be$.
In fact, not much ordinal theory is presumed in this article; we almost exclusively work with ordinals less than $\varphi_20$, which can be written in terms of addition and the functions $\xi\mapsto\omega^\xi$, $\xi\mapsto \ve_\xi$.
It will be convenient for us to adopt the convention that $\ve_{-1} = 0$.
We assume familiarity with ordinal addition and multiplication.
We also use partial subtraction: recall that if $\alpha>\beta$ there is a unique ordinal $\eta$, such that $\beta + \eta =\alpha$.
We will denote this unque $\eta$ by $-\beta+\alpha$
For more details, we refer the reader to standard texts such as \cite{Pohlers:2009:PTBook,Schutte1977}.

\begin{definition}
Let $\Lambda$ be an ordinal. A {\em system of fundamental sequences on $\Lambda$} is a function $\cdot [\cdot] \colon \Lambda \times \mathbb N \to \Lambda$ such that $  \alpha [n] \leq \alpha$ with equality holding if and only if $\alpha = 0$, and $\alpha[n] \leq \alpha[m]$ whenever $n\leq m$.
The system of fundamental sequences is {\em convergent} if $\lambda  = \lim_{n\to \infty} \lambda[n]$ whenever $
\lambda $ is a limit, and has the {\em Bachmann property} if whenever $\alpha[n] < \beta < \alpha$, it follows that $\alpha[n] \leq \beta[1]$.
\end{definition}

Define $\fsi \la {n} = \la [0][1] \ldots [n]$.
It is clear that if $\Lambda$ is an ordinal then for every $\alpha<\Lambda $ and every $m$ there is $n$ such that $  \fsi \al {n} = 0$, but this fact is not always provable in weak theories.
The ordinal $\varphi_20$ enjoys a natural system of fundamental sequences satisfying  the Bachmann property~\cite{Schmidt77}.

\begin{definition}
Let $\om_0(\al):=\al$ and $\om_{k+1}(\al)=\om^{\om_k(\al)}$.
The standard fundamental sequences for ordinals up to $\varphi_20$ are defined as follows.
\begin{enumerate}

\item If $\al= \omega^\beta+\gamma$ with $0<\gamma<\alpha$, then
$\al[k]:= \omega^\beta+\gamma[k]$.

\item If $\al=\om^\be>\be$, then we set $\al[k]:= 0$ if $\be=0$, $\al[k]:=\om^\ga\cdot k$ if $\be=\ga+1$, and $\al[k]:=\om^{\be[k]}$ if $\be\in {\rm Lim}$.

\item If $\al=\varepsilon_\be>\be$, then 
$\al[k]:=\om_k(1)$ if $\be=0$, $\al[k]:=\om_k(\varepsilon_\ga+1)$ if $\be=\ga+1$, and $\al[k]:=\varepsilon_{\be[k]}$ if $\be\in{\rm Lim}$.

\item If $\al=\varphi_20$ then $\al[0]=1$ and $\al[n+1] = \ve_{\al[n]}$. 

\end{enumerate}
\end{definition}

We will use these fundamental sequences to establish new independence results by appealing to proof-theoretic ordinals.
Such ordinals are defined in terms of a given ordinal notation system with fundamental sequences, and in this case, we use our notation system for $\varphi_20$.
Each of the theories $T$ we are interested in can be assigned a proof-theoretic ordinal $\|T\| \leq \varphi_20$, which gives us a wealth of information about what is (un)provable in $T$.
First, it characterizes the amount of transfinite induction available.
Recall that transfinite induction for an ordinal $\al$ is the scheme
\[{\rm TI}(\al) := \forall \xi < \bar\al \Big (\big (\forall \zeta<\xi \ \varphi(\zeta) \big) \to \varphi (\xi) \Big )\rightarrow \forall \xi<\bar\al \ \varphi (\xi),\]
where for our purposes $\varphi$ is a formula in the language of Peano arithmetic and $\bar\al$ is the numeral coding $\al$.
Then, $\|T\|$ can be characterized as the least ordinal $\xi$ such that $T\not\vdash {\rm TI}(\xi)$.

The ordinal $\|T\|$ can also be used to bound the provably total computable functions in $T$.\footnote{It should be remarked that proof-theoretic ordinals defined in terms of transfinite induction or in terms of computable functions are not necessarily equivalent, but they tend to coincide for `natural' theories, including those considered in this text.}
Given $\xi\leq\varphi_2 0 $, let $H(\xi)$ be the least $n$ so that $\fsi\xi n =0$, and for $n\in\mathbb N$, let $F_\xi(n) = H ( \xi [n])$.
Then, the proof-theoretic ordinal of $T$ is also the least $\xi\leq \varphi_2(0)$ such that $T\not\vdash \forall \zeta<\xi \ \forall n \  \exists m \ ( m = F_\xi(n) )$, if it exists (otherwise, we may define it to be $\infty$, or just leave it undefined).

Say that a partial function $f\colon \mathbb N \to \mathbb N$ is {\em computable} if there is a $\Sigma_1$ formula $\varphi_f (x,y)$ in the language of first order arithmetic (with no other free variables) such that for every $m,n$, $f(m) = n$ if and only if $\varphi_f (m,n)$ holds.
The function $f$ is {\em provably total} in a theory $T$ if $T\vdash \forall x \exists y \varphi _f (x,y) $ (more precisely, $f$ is provably total if there is at least one such choice of $ \varphi _f $).

\begin{theorem}\label{theoIndep}
Define
\begin{itemize}

\item $\|{\sf ACA} _0\| =\ve_0$ \cite{Pohlers:2009:PTBook},

\item $\|{\sf ACA}' _0\| =\ve_\om$ \cite{Afshari12}, and

\item $\|{\sf ACA}^+ _0\| = \varphi_2(0)$ \cite{Leigh10}.

\end{itemize}
Then, for $T\in \{ {\sf ACA} _0 , {\sf ACA}' _0,{\sf ACA}^+ _0\}$ and $\al<\varphi_2(0)$, the following are equivalent~\cite{BCW}:
\begin{enumerate}

\item $\al<\|T\|$.

\item $T\vdash {\rm TI}( \al)$.

\item $F_\al$ is provably total in $T$.

\item There exists a provably total computable function $f\colon\mathbb N\to\mathbb N$ in $T$ such that $f(n) > F_\al(n)$ for all $n$.

\end{enumerate}
\end{theorem}

The Bachmann property will be useful in transferring unprovability results stemming from proof-theoretic ordinals to the setting of Goodstein processes, in view of the following.

\begin{proposition}\label{propMajorize}
Let $\Lambda$ be an ordinal with a system of fundamental sequences satisfying the Bachmann property, and let $(\xi_n)_{n\in \mathbb N}$ be a sequence of elements of $\Lambda$ such that, for all $n$, $\xi_n[n+1] \leq \xi_{n+1} \leq \xi_n$.
Then, for all $n$, $\xi_n \geq \fsi{\xi_0}n$.
\end{proposition}

\proof
           Let $\preceq_k$ be the reflexive transitive closure of $\{(\al[k],\al):\al<\varphi_2(0)\}$.
           We need a few properties of these orderings.
           Clearly, if $\al \preceq_k \be$, then $\al \leq \be$.
                     It can be checked by a simple induction and the Bachmann property that, if $\alpha[n] \leq \beta < \alpha$, then $\alpha[n] \preceq_1 \beta$.
Moreover, $\preceq_k$ is monotone in the sense that if
          $\al\preceq_k \be$, then $\al\preceq_{k+1}\be$, and if $\alpha\preceq _k \beta$, then $\alpha[k] \preceq_k \beta[k]$ (see, e.g., \cite{Schmidt77} for details).

We claim that for all $n$, $\xi_n \succeq _n \xi_0[1]\ldots[n]$, from which the desired inequality immediately follows.
For the base case, we use the fact that $\succeq_0$ is transitive by definition.
For the successor, note that the induction hypthesis yields $\xi_0[1]\ldots[n] \preceq_n \xi_n $, hence $\xi_0[1]\ldots[n+1] \preceq_{n+1} \xi_n[n+1] $.
Then, consider three cases.

\begin{enumerate}[label*={\sc Case \arabic*},wide, labelwidth=!, labelindent=0pt]

\item ($\xi_{n+1} = \xi_n$).
By transitivity and monotonicity, $\xi_0[1]\ldots [n+1] \preceq_{n+1} \xi_0[1]\ldots [n ] \preceq _n \xi_n =  \xi_{n+1}$
yields $\xi_0[1]\ldots [n+1] \preceq_{n+1}  \xi_{n+1}$.

\item ($\xi_{n+1} = \xi_n[n+1]$).
Then, $\xi_0[1]\ldots [n+1] \preceq_{n+1} \xi_n[n+1] = \xi_{n+1}$.

\item ($\xi_n[n+1] < \xi_{n+1} < \xi_n$).
The Bachmann property yields $\xi_n[n+1] \preceq_1 \xi_{n+1}$, and since $\xi_0[1]\ldots[n+1] \preceq_{n+1} \xi_n[n+1] $, monotinicity and transitivity yield $\xi_0[1]\ldots[n+1] \preceq_{n+1} \xi_{n+1} $.\qedhere

\end{enumerate}
\endproof

As an immediate corollary we obtain that for such a sequence of ordinals $(\xi_n)_{n\in\mathbb N}$, if $m$ is least so that $\xi_m = 0$, then $m\geq F_{\xi_0}(n)$.
So, our strategy will be to show that the respective Goodstein process for $T$ grows at least as quickly as $F_{\|T\|}$, from which we obtain that the theorem is unprovable in $T$.

\section{Goodstein sequences for ${\sf ACA}_0$}

Each of the three notions of nested normal form we have discussed (nested, index-nested, and argument-nested) leads to a Goodstein principle of differing proof-theoretical strength.
The key here is that each notion of normal form leads to a different base change operation, and each base change operation leads to a faster- or slower-terminating Goodstein process.
In this section we study argument-nested normal forms, and show that they lead to a Goodstein principle of the same proof-theoretic strength of the original.

\begin{definition}\label{defBCHA}
For $2\leq k \leq \ell$ and $m\in \mathbb N$, define $m\bch{k}{\ell}$ by:
\begin{enumerate}
\item $0\bch{k}{\ell}:=0$
\item $m\bch{k}{\ell}:=A_{a\bch{k}{\ell}}(\ell,b) +c\bch{k}{\ell}$ if $m\gnf A_a(k,b) + c$.
\end{enumerate}
We will write $m\bchk$ instead of $m\bch k{k+1}$.
\end{definition}

It is not hard to check that if $m\gnfm A_a(k,b)\cdot d+e$, then
\[m\bch{k}{\ell} = A_{a\bch{k}{\ell}} (\ell,b) \cdot d + e\bch{k}{\ell}.\]
We then define a new Goodstein process based on this new base change operator.

\begin{definition}
Let $\ell<\om$. 
Put $\IG_0(\ell):=\ell.$
Assume recursively that $\IG_k(\ell)$ is defined and $\IG_k(\ell)>0$.
Then, $\IG_{k+1}(\ell)=\IG_k(\ell)\bch{k+2}{k+3}-1$. If $\IG_k(\ell)=0$, then $\IG_{k+1}(\ell):=0$.
\end{definition}

We will show that for every $\ell$ there is $i$ with $\IG_i(\ell) = 0$.
In order to prove this, we first establish some natural properties of the base-change operation.

\begin{lemma}\label{lemUnMonA}
Fix $k\geq 2$ and let $m,n\in\mathbb N$.
\begin{enumerate}
\item $m\leq m\bchk$.
\item If $m<n$, then $m\bchk<n\bchk$.
\end{enumerate}
\end{lemma}

\proof
Write $A_ab$ for $A_a(k,b)$ and $B_ab$ for $A_a(k+1,b)$.
The first assertion is proved by induction on $m$. It clearly holds for $m=0$.
If $m\gnf A_ab +c$ then the induction hypothesis yields
$m= A_ab +c \leq B_{ a\bchk}(b) +c\bchk=m\bchk.$

For the second we proceed by induction on $n$ with a secondary induction on $m$.
The assertion is clear if $m=0$.
Let $m\gnf A_a b  +c$
and $n\gnf A_{{a'}} {b'}  +c'$.
First we establish some useful inequalities.
Note that $A_a b \leq m < A_{a+1} 0  = A_{a}^k0= A_aA_a^{k-1} 0 $. So $b < A_a^{k-1}0 \leq B_{a\bchk}^{k-1}0$.
Moreover $c<A_a(b+1) $ yields
\[c\bchk<B_{a\bchk}(b+1) \leq B_{a\bchk}  A_{a\bchk}^{k-1}0  \leq B_{a\bchk}^{k}0.\]
Now, consider several cases.

\begin{enumerate}[label*={\sc Case \arabic*},wide, labelwidth=!, labelindent=0pt]
\item ($a<{a'}$).
By the induction hypothesis $a <{a'}$ yields $(a+1)\bchk \leq {a'}\bchk$.
This yields
\begin{align*}
m\bchk & = B_{a\bchk }b +c\bchk
<  B_{a\bchk}  B_{a\bchk }^{k-1}0 + B_{a\bchk}^{k}0\\
&= 2 B_{a\bchk}^{k}0 \leq B_{a\bchk}^{k+1}0 = B_{a\bchk+1} 0
\leq B_{(a+1)\bchk } 0 \\
&\leq B_{{a' }\bchk }(b') +c'\bchk  = n\bchk,
\end{align*}
where the second equality uses Lemma \ref{lemmUseful}.\ref{itMult}.

\smallskip \item ($a={a'}$ and $b<{b'}$).
We consider two sub-cases:
\smallskip

\begin{enumerate}[label*= .\arabic*,wide, labelwidth=!, labelindent=0pt]
\item ($A_a(b+1)<n$).
Since $m < A_a(b+1)$ and $b+1\leq {b'}$ then by the induction hypothesis
$m\bchk<B_{a\bchk} (b+1)  \leq B_{{a'}\bchk}{b'} \leq n\bchk$.

\smallskip \item ($A_a( b+1)=n$). Write $m = A_ab \cdot d + e$ and we consider two further sub-cases depending on the value of $n$.

\smallskip

\begin{enumerate}[label*= .\arabic*,wide, labelwidth=!, labelindent=0pt]
\item
($a=0$).
In this case we may write $m\gnfm k^b\cdot d+e$, with $d<k$ and $e<k^b$, and $n$ has $k$-normal form $k^{b+1}$.
The induction hypothesis applied to $e<k^b$ yields $e \bchk < (k+1)^b$.
We then have that $m\bchk=(k+1)^b\cdot d +e\bchk < (k+1)^b\cdot k + (k+1)^b = (k+1)^{(b+1)} = n\bchk$.
\smallskip

\item($a>0$). Write $m\gnfm A_ab\cdot d+ e$.
Note that $d<A_a(b+1)$ and $e<A_ab$, which by the induction hypothesis yields $e\bchk <B_{a\bchk}(b ) $.
Let $r= B_{a\bchk-1}^kB_{a\bchk}(b )$, so that $B_{a\bchk-1} r = B_{a\bchk}(b+1)$.
Then,
\begin{align*}
m\bchk 
&=B_{a\bchk}  (b)\cdot d  +e\bchk
\leq  B_{a\bchk}(  b ) \cdot A_a(b+1) +B_{a\bchk}( b )\\
& < r^2+r
\leq B_{a\bchk-1} r = B_{a\bchk}( b + 1 ) \leq B_{a\bchk}( (b+1)  )= n\bchk,
\end{align*}
where the second inequality follows by
\[A_a(b+1) = A_{a-1}^k A_ab \leq B_{a-1}^kB_{a\bchk}(b ) = r \]
and the third inequality uses Lemma \ref{lemmUseful}.\ref{itSquare}.

\end{enumerate}
\end{enumerate}

\item ($a={a'}$ and ${b'}=b$). Since $m<n$ then we must have that $c<c'$, and
\[
m\bchk
=B_{a\bchk}b+c\bchk\\
<B_{a\bchk}b+c'\bchk
= n\bchk .
\qedhere\]

\end{enumerate}
\endproof

Thus, the base-change operation is monotone.
Next we see that it also preserves normal forms.

\color{black}
\begin{lemma}\label{prA}
Fix $k\geq 2$.
If $m=A_a(k,b) +c$ is in $k$-normal form, then $m \bchk=A_{a\bchk}(k+1,b) +c\bchk$ is in $k+1$ normal form.
\end{lemma}

\proof
Write $A_xy$ for $A_x(k,y)$ and $B_xy$ for $A_x(k+1,y)$.
Let $m\gnf A_a b+c$.
We have that  $m<A_{a+1} 0$,  $m<A_a( b+1)$, and $c<A_ab$.
So, $A_ab < A_{a+1}0 < A_a^k0 = A_a A_a^{k-1}0$.
Hence $b< A_a^{k-1}0 < B_{a\bchk}^{k-1}0 ; $
Since $A_a b $ is in $k$-normal form, Lemma \ref{lemUnMonA} yields
$c\bchk < B_{a\bchk}b < B_{a\bchk}B_{a\bchk}^{k-1}0 $.
So 
\begin{align*}
m\bchk&
<B_{a\bchk } B_{a\bchk}^{k-1}0 + B_{a\bchk}B_{a\bchk}^{k-1}0\\
&= 2  B_{a\bchk}^{k}0 \leq B_{a\bchk}^{k+1}0 = B_{a\bchk +1}0\leq B_{(a+1)\bchk }0.
\end{align*}

Now we check that $ m\bchk < B_{a\bchk}( b + 1)$.

If $a=0$, then $m \gnf A_0b +c \gnfm k^b\cdot d+e$ for some $d<k$ and $e<k^b$.
Note that $k^b\cdot d'+e$ is in extended normal form if $1 \leq d'<d$, from which it readily follows that
\begin{align*}
m\bchk & = (k+1)^b\cdot d+e\bchk < (k+1)^b\cdot k+(k+1)^b \\
& =(k+1)^{b+1} = B_0(b+1).
\end{align*}

In the remaining case, write $m\gnfm A_ab \cdot d + e$, with $e<A_ab$, so that $d<A_a(b+1)$.
By Lemma \ref{lemmUseful}.\ref{itSquare}, $B_{a\bchk}(b+1) \geq B_{a\bchk}B_{a\bchk}b\geq (B_{a\bchk}b)^2 +B_{a\bchk}b $, so
\begin{align*}
m\bchk
&=B_{a\bchk}b\cdot d+e\bchk
<B_{a\bchk} b  \cdot A_{a\bchk}(b+1) + B _{a\bchk}b \\
& \leq (B_{a\bchk} b)^2+B _{a\bchk}b
\leq B_{a\bchk}( b+1  ).\\
\end{align*}
So, $A_{a\bchk}(k+1,b) +c\bchk$ is in $k+1$-normal form.
\endproof

\begin{definition}
For $k\geq 2$, define $\zk\cdot \colon \mathbb N \to \ve_0$ as follows:
\begin{enumerate}
\item $\zk0:=0$.
\item If $m\gnf A_ab+c$, then $\zk m:=\om^{\om \cdot (\zk a) + b} +\zk c$.
\end{enumerate}
\end{definition}

The function $\zk {\cdot}$ is the {\em base $k $ ordinal assignment,} and provides a monotone mapping of the natural numbers into the ordinals.

\begin{lemma}\label{lemPkMonA}
If $m<n<\om$ then $\zk m<\zk n$.
\end{lemma}

\proof
Proof by induction on $n$ with subsidiary induction on $m$.
The assertion is clear if $m = 0$.
Let $m \gnf A_ab+c$
and $n \gnf A_{{a'}} {b'} +c'$.

\begin{enumerate}[label*={\sc Case \arabic*},wide, labelwidth=!, labelindent=0pt]
\item ($a<{a'}$). 
By the induction hypothesis $\zk a< \zk {{a'}}$, thus $\om^{\om \cdot \zk a}<\om^{\om \cdot \zk {{a'}}}$.
We have $c<m<A_{a+1}0\leq A_{{a'}} 0 \leq n$, and the induction hypothesis yields
$\zk c < \zk{(A_{{a'}}0)} = \om^{\om \cdot \zk {{a'}}}$.\\
Since $ \om^{\om \cdot \zk a +b} <  \om^{\om \cdot \zk {{a'}}}$ then
$\zk m =  \om^{\om \cdot \zk a + b } + \zk c < \om^{\om \cdot \zk {{a'}}} \leq \zk n$.
\smallskip

\item ($a={a'}$). We consider several sub-cases.
\smallskip

\begin{enumerate}[label*= .\arabic*,wide, labelwidth=!, labelindent=0pt]
\item ($b<{b'}$). 
By the induction hypothesis $\om^{\om \cdot \zk a +b}< \om^{\om \cdot \zk a + {b'}}$.
We have  $c<A_a b $, and the subsidiary induction hypothesis yields
\[\zk c < \zk {(A_ab)} = \om^{\om \cdot \zk a +b} < \om^{\om \cdot \zk a + {b'}}.\]
So,
$\zk m = \om^{\om \cdot \zk a +b}  + \zk c < \om^{\om \cdot \zk a + {b'}} \leq \zk n.$
\smallskip

\item ($b={b'}$).
Then $c<c'$ and by the induction hypothesis $\zk c< \zk {c'}$.
Hence 
$\zk m = \om^{\om \cdot \zk a +b}  + \zk c < \om^{\om \cdot \zk {{a'}} +{b'}}   + \zk {c'} = \zk n.$

\end{enumerate}
\end{enumerate}
\endproof

The general idea of the proof of termination is that the functions $\zk\cdot$ can be used to assign a decreasing sequence of ordinals to any Goodstein sequence.
Since no decreasing sequence of ordinals can be infinite, the sequence must terminate at some point.
For this to work, we must show that the ordinal assignment is monotone, in the following sense.

\begin{lemma} 
For any $m\in  \mathbb N$ and $k\geq 2$, $\zskpe { m\bch k{k+1} }=\zk m$.
\end{lemma}

\proof
Proof by induction on $m$.
The assertion is clear for $m=0.$
Assume $m\gnf A_a(k,b) +c$.
Then, $m\bchk \gnfp A_{a\bchk}(k+1,b) +c\bchk$, and the induction hypothesis yields

\begin{eqnarray*}
\zskpe { m\bchk }&=&\zskpe {(A_{a\bchk}(k+1,b) +c\bchk)}\\
&=&\om^{\om \cdot \zskpe { a\bchk } +b } +\zskpe {c\bchk}\\
&=&\om^{\om \cdot \zk a +b} +\zk c =\zk m.
\end{eqnarray*}
\endproof

\begin{lemma}\label{majA}
Given $k,m<\omega$ with $k\geq 2$, $\zskpe {(m\bchk-1)} \geq  \zk m [k]$.
\end{lemma}

\proof
We prove the claim by induction on $c$.
Write $A_xy$ for $A_x(k,y)$ and $B_xy$ for $A_x(k+1,y)$.
Let $m\gnf A_ab +c.$ 
\smallskip

\begin{enumerate}[label*={\sc Case \arabic*},wide, labelwidth=!, labelindent=0pt]

\item ($c>0$).
Then the induction hypothesis and Lemma \ref{lemPkMonA} yields
\begin{align*}
\zskpe {(m\bchk-1)}
&= 
\zskpe {(  B_{a\bchk}b +c \bchk -1 )}\\
&=  \om^{\om \cdot \zskpe {a\bchk} +b } +(\zskpe{c\bchk -1})\\
&\geq \om^{\om \cdot \zk a + b} + \zk c [k]\\
&=   \zk{( A_ab + c)} [k]
= \zk m[k].
\end{align*}

\item ($c=0$).
We consider several sub-cases.
\smallskip

\begin{enumerate}[label*= .\arabic*,wide, labelwidth=!, labelindent=0pt]
\item ($a>0$ and $b>0$)
We have $m=A_ab$. Then by the induction hypothesis and Lemma \ref{lemPkMonA},
\begin{align*}
\zskpe {(m\bchk-1)}
&= \zskpe{(B_{a\bchk}b-1)}\\
&\geq \zskpe{B_{a\bchk}( b-1)} \\
&= \om^{\om \cdot \zskpe{a\bchk} + (b-1)}\\
&\geq \om^{\om \cdot \zk a + (b-1)}\\
&=  \om^{\om \cdot \zk a +b[k]} 
=  \zk m [k].
\end{align*}
\smallskip

\item  ($a>0$ and $b=0$).
In this case $m=A_a 0$, and the induction hypothesis yields
\begin{align*}
\zskpe{(m\bchk-1)}&  =  \zskpe{(B_a0-1)}\\
&= \zskpe{(B_{a\bchk-1}^{k+1}0-1)}\\
& \geq  \zskpe{(B_{a\bchk-1}^{k}0)}\\
&= \om^{\om \cdot \zskpe {(a\bchk -1)} + B_{a\bchk -1}^{k-1}0}\\
&\geq   \om^{\om \cdot \zk  a [k] + B_{a\bchk -1}^{k-1}0}\\
&\geq \om^{\om \cdot \zk (a)[k]}
=   \zk{A_a0}[k].
\end{align*}
\smallskip

\item  ($a=0$ and $b>0$).
Then $m=A_0b$. 
\begin{align*}
\zskpe{(m\bchk-1)}
&= \zskpe {(B_0b-1)}\\
&\geq \zskpe {( B_0( b-1)} \cdot k)\\
&= \om^{b-1} \cdot k\\
&=\om ^{b[k]}
= (\zk c)[k].
\end{align*}
since $B_0( b-1) \cdot k$ is in $k+1$ normal form.
\smallskip

\item ($a=0$ and $b=0$)
In this case
$\zk m [k] = \zk{B_00} [k] = (\omega^{\omega\cdot 0+0})[k] = 1[k] = 0$, so the lemma follows.

\end{enumerate}
\end{enumerate}
\endproof

The proof of termination follows by observing that for every $\ell$,
\[ \io_k(\ell ):= \zk{\IG_k(\ell)} \]

is a decreasing sequence of ordinals, hence must be finite. Below, we make this precise.
It is well-known that the so-called slow-growing hierarchy at level $\varphi_\omega 0$ matches up with the Ackermann function, so one might expect that the corresponding Goodstein process can be proved terminating in ${\sf PA}+{\rm TI}(\varphi_\omega0)$.
This is true but, somewhat surprisingly, much less is needed here. We can lower $\varphi_\omega0$ to $\varepsilon_\om=\varepsilon_0$.

\begin{theorem}\label{ThmTermA}
For all $\ell<\om$, there exists a $k<\om$ such that $\IG_k(\ell)=0.$ This is provable in ${\sf PA}+{\rm TI}(\ve_0)$.
\end{theorem} 

\proof
If $\IG_k(\ell) \bch \omega k>0$, then, by the previous lemmata,
\begin{eqnarray*}
\io_{k+1}(\ell )&=&\IG_{k+1}(\ell) \bch {k+3}\omega 
= (\IG_{k}(\ell)\bch{k+2}{k+3}-1)\bch {k+3}\omega \\
&<& \IG_{k}(\ell)\bch{k+2}{k+3} \bch {k+3}\omega
=\IG_{k}(\ell)  \bch {k+2}\omega=\io_k(\ell).
\end{eqnarray*}
Since $(\io_k(\ell))_{k<\omega}$ cannot be an infinite decreasing sequence of ordinals, there must be some $k$ with $\io_k(\ell) = 0$, yielding $\IG_k(\ell) = 0$.
\endproof

\begin{theorem}\label{TheoACA}
${\sf ACA}_0\not\vdash \forall \ell \exists k\ \IG_k(\ell )=0.$
\end{theorem}

\proof
Recall that $\io_k(A_\ell (2,0) )=\IG_k(A_\ell (2,0))\bch {k+2}\omega$. 
Let $\ell_n$ be given by $\ell_0 = 0$ and $\ell_{n+1} = A_{\ell_n} (2,0)$.
Then, we have that $\io_0 (\ell_{n+1}) = \omega_n(1)$.
Lemma \ref{majA} and Lemma \ref{lemPkMonA} yield $\io_k(A_{\ell_{n+1}} (2,0) )[k+1] \leq \io_{k+1}(A_{\ell_{n+1}} (2,0) ) < \io_k (A_{\ell_{n+1}} (2,0)  )$, hence Proposition \ref{propMajorize} yields $ \io_k (\ell_{n+1} ) \geq \fsi{\io_0 (A_\ell (2,0)) }k.$
So the least $k$ such that $\IG_k(\ell_{n+1})=0$ is at least as big as $H(\omega_n ) = F_{\ve_0}(n)$.
By Theorem \ref{theoIndep}, $\|{\sf ACA}_0\| = \ve_0$, and thus ${\sf ACA}_0$ does not prove that this $k$ is always defined as a function of $\ell_{n+1}$, so that ${\sf ACA}_0\not \vdash \forall \ell \exists k \ \IG_k(A_\ell (2,0)))=0.$
\endproof

\section{Goodstein sequences for ${\sf ACA}_0'$}\label{secPrime}

The base change operator of the previous section corresponds to the notion that $A_ab+c$ is written in such a way that $a$ and $c$ are also expressed using the Ackermann function, but $b$ is treated as a parameter.
But what if we instead treat $a$ as a parameter and express $b$ using the Ackermann function? One may expect this to lead to a {\em weaker} Goodstein result, since the function $a\mapsto A_ab$ grows much faster than $b\mapsto A_ab$.
Surprisingly, this is not the case, as the new Goodstein process is independent of ${\sf ACA}'_0$.

\begin{definition}
Given $k\geq 2$, $m\in\mathbb N$ and $2\leq k<\ell$, define $m\bchb k\ell$ by:
\begin{enumerate}
\item $0\bchb k\ell:=0$.
\item $m\bchb k\ell :=A_a \left (\ell,b\bchb k\ell\right ) +c\bchb k\ell$ if $m\gnf A_a(k,b) +c$.
\end{enumerate}
We write $\bchkb$ instead of $\bchb k{k+1}$.
\end{definition}

The idea is that now, $a$ is being treated as a constant, while $b$ should be regarded as being iteratively represented using the Ackermann function.
As before, if $m\gnfm A_a(k,b)\cdot d+e$, then
\[m\bch{k}{\ell} = A_{a} (\ell,b\bchb{k}{\ell}) \cdot d + e\bchb{k}{\ell}.\]
We then define a new Goodstein process based on this new base change operator.
This new base change operator leads to a second Goodstein process.

\begin{definition}
Let $\ell<\om$. 
Put $\BG_0(\ell):=\ell.$
Assume recursively that $\BG_k(\ell)$ is defined and $\BG_k(\ell)>0$.
Then $\BG_{k+1}(\ell)=\BG_k(\ell)\bchb{k+2}{k+3}-1$. If $\BG_k(\ell)=0$, then $\BG_{k+1}(\ell):=0$.
\end{definition}

We will prove monotonicitiy using a similar method as in the previous section.

\begin{lemma}\label{lemUnMon}
Fix $k\geq 2$ and let $m,n\in\mathbb N$. Then:
\item If $m<n$, then $m\bchkb<n\bchkb$.
\end{lemma}

\proof
Write $A_xy$ for $A_x(k,y)$ and $B_xy$ for $A_x(k+1,y)$.

The proof is by induction on $n$ with a subsidiary induction on $m$.
The assertion is clear if $m=0$.
Let $m\gnf A_ab +c$
and $n\gnf A_{{a'}} {b'}  +c'$.
We distinguish cases according to the position of $a$ relative to ${a'}$, the position of $b$ relative to ${b'}$, etc.
Note that from the choice of $a,{a'}$ we must have $a\leq a '$.
\begin{enumerate}[label*={\sc Case \arabic*},wide, labelwidth=!, labelindent=0pt]
\item ($a<{a'}$).
Write $m\gnfm A_ab\cdot d+e$.
We have $A_ab\leq m < A_{a+1}0 = A_a A_a^{k-1}0$.
For $\ell <k$ we have that $A_a^{\ell }0$ is in $k$-normal form by Lemma \ref{NF}.
Thus the induction hypothesis yields $b \bchkb < B_a^{k-1}0$.
The number
$A_ab$ is in $k$-normal form and so the induction hypothesis applied to $e<A_ab$ yields
$e\bchkb<B_a  (b\bchkb) \leq B_a^{k }0$.
Moreover, $d< m < A_{a+1}0 $.
This yields
\begin{align*}
m\bchkb
& = B_{a}(b\bchkb  )\cdot d + e\bchkb \\
& \leq  B_{a} B_a^{k-1}0 \cdot A_a^{k}0  + B_aB_a^{k-1}0 \\
& \leq (B_a^{k }0)^2 + B_a^{k }0 < B_a^{k +1}0 = B_{a+1}0\\
&\leq B_{a'}(b\bchkb) + c'\bchkb =n\bchkb,
\end{align*}
where the second inequality follows from
\[A_{a+1}0 = A_a^{k}0 \leq B_a^{k} 0 \]
and the second-to-last from Lemma \ref{lemmUseful}.\ref{itSquare}.

\item ($a={a'}$). Note that in this case $b\leq {b'}$. Consider the following sub-cases.

\begin{enumerate}[label*= .\arabic*,wide, labelwidth=!, labelindent=0pt]

\smallskip \item ($a={a'}$ and $b<{b'}$).
Consider two sub-cases.
\smallskip

\begin{enumerate}[label*= .\arabic*,wide, labelwidth=!, labelindent=0pt]
\item ($A_a( b+1)<n$).
Since $n$ is in $k$-normal form and $b+1\leq {b'}$ we see that $A_a( b+1)$ is in $k$-normal form by Lemma \ref{NF}.
Then, the induction hypothesis yields
$m\bchkb < B_a ((b+1)\bchkb) \leq B_a( {b'}\bchkb) \leq n\bchkb$.

\smallskip \item ($A_a( b+1)=n$).
We know that $m=A_ab +c<A_a( b+1)=n$.
Consider two further sub-cases.
\smallskip

\begin{enumerate}[label*= .\arabic*,wide, labelwidth=!, labelindent=0pt]
\item
($a=0$). This means that $m \gnfm k^b\cdot d +e < k^{b+1}=n$ for some $d<k$ and $e<k^b$,
where $n$ has $k$-normal form $k^{b+1}$.
The induction hypothesis yields $b\bchkb<(b+1)\bchkb$ and $e\bchkb<(k+1)^{b\bchkb}$.
We then have that
\begin{align*}
m\bchkb&=(k+1)^{b\bchkb} \cdot d +e\bchkb <(k+1)^{b\bchkb  } \cdot k + (k+1)^{b\bchkb  }\\
&  \leq (k+1)^{b\bchkb +1 }\leq  (k+1)^{(b+1)\bchkb  } = n\bchkb.
\end{align*}

\item($a>0$).
Write $m\gnfm A_ab\cdot d+ e$.
Note that $d<A_a(b+1)$ and $e<A_ab$, which by the induction hypothesis yields $e\bchkb <B_a(b\bchkb ) $.
Let $r= B_{a-1}^kB_a(b\bchkb )$, so that $B_{a-1} r = B_a(b\bchkb +1)$.
Then,
\begin{align*}
m\bchkb 
&=B_{a}  (b\bchkb)\cdot d  +e\bchkb
\leq  B_{a}(  b\bchkb ) \cdot A_a(b+1) +B_a( b\bchkb )\\
& < r^2+r
\leq B_{a-1} r = B_a( b\bchkb + 1 ) \leq B_a( (b+1)\bchkb  )= n\bchkb,
\end{align*}
where the second inequality follows by
\[A_a(b+1) = A_{a-1}^k A_ab \leq B_{a-1}^kB_a(b\bchkb ) = r \]
and the third inequality uses Lemma \ref{lemmUseful}.\ref{itSquare}.
\end{enumerate}

\end{enumerate}

\smallskip \item (${b'}=b$). Then $c<c'$ and the induction hypothesis yields
\begin{align*}
m\bchkb
&=B_{a}( b\bchkb  ) +c\bchkb\\
&<B_{a}( b\bchkb  ) +c'\bchkb
= n\bchkb . &\qedhere
\end{align*}
\end{enumerate}
\end{enumerate}
\endproof

Thus, the base-change operation is monotone.
Next we see that it also preserves normal forms.

\begin{lemma}\label{pr}
If  $m \gnf A_a(k,b) +c$, then
 $m\bchkb \gnfp A_a(k+1,b\bchkb) +c\bchkb$.
\end{lemma}

\proof
Write $A_xy$ for $A_x(k,y)$ and $B_xy$ for $A_x(k+1,y)$.
Assume that $m\gnf A_a(k,b) +c$.
Then, $m <A_{a+1} 0$, 
$m<A_a(b+1)$, and $c<A_ab$.
Clearly, $ B_a0 \leq m\bchkb $.
By Lemma \ref{NF}, $A_{a+1}0$ is in $k$-normal form, so that by Lemma \ref{lemUnMon}, 
$c<A_{a+1} 0 $ yields $c\bchkb<B_{a+1} 0 $.
Since $A_ab$ is in $k$-normal form, Lemma \ref{lemUnMon} yields
$c\bchkb < B_a( b\bchkb)$.
It remains to check that we also have $ m\bchkb < B_a( b\bchkb + 1)$.

If $a=0$, then write $m \gnfm A_0b\cdot d +e$ for some $d<k$ and $e< k^b$.
Then, $m\bchkb=(k+1)^{b\bchkb}\cdot d+e\bchkb<(k+1)^{b\bchkb+1}$
and thus $m\bchkb\gnf (k+1)^{b\bchkb} +c\bchkb.$
In the remaining case, we have for $a>0$ that
\begin{align*}
m\bchkb  &
=B_{a}(  b\bchkb  ) +c\bchkb 
<B_{a}( b\bchkb  ) +B_{a}( b\bchkb  ) \\
& = 2 B_{a}( b\bchkb  )
 \leq  B_{a}(  b\bchkb+1  ),
\end{align*}
with the last equality following from Lemma \ref{lemmUseful}.\ref{itMult}.
So $B_a( b\bchkb) +c\bchkb$ is in $k+1$-normal form.
\endproof

The Goodstein process arising from this base-change operator is also terminating.
In order to prove this, we must assign ordinals to natural numbers, in such a way that the process gives rise to a decreasing (hence finite) sequence.
For each $k$, we define a function $\pk\cdot \colon \mathbb N \to \Lambda$, where $\Lambda$ is a suitable ordinal, in such a way that $\pk m$ is computed from the $k$-normal form of $m$.
Unnested Ackermannian normal forms correspond to ordinals below $\Lambda=\ve_\omega$, as the following map shows.
Below, recall that as per our convention, $\ve_{-1} = 0$.
\begin{definition}
For $k\geq 2$, define $\pk \cdot \colon \mathbb N \to \ve_\omega$ as follows:
\begin{enumerate}
\item $\pk 0:=0$.
\item $\pk m:=\om^{\ve_{a-1}+\pk b} +\pk c$ if $m\gnf A_a(k,b) + c$.
\end{enumerate}
\end{definition}

The ordinal assignment in this case is once again monotone.

\begin{lemma}\label{lemPkMon}
If $m<n<\om$ then $\pk m<\pk n$.
\end{lemma}

\proof
Proof by induction on $n$ with subsidiary induction on $m$.
The assertion is clear if $m=0$.
Write $A_xy$ for $A_x(k,y)$ and let $m\gnf A_a b +c$
and $n\gnf A_{{a'}} {b'}  +c'$.
We distinguish cases according to the position of $a$ relative to ${a'}$, the position of $b$ relative to ${b'}$, etc.
Note that $a\leq {a'} $.

\begin{enumerate}[label*={\sc Case \arabic*},wide, labelwidth=!, labelindent=0pt]
\item ($a<{a'}$). We have $c<m<A_{a+1} 0 \leq A_{{a'}} 0 $ and, since $A_{{a'}}0 \leq n$, the induction hypothesis yields
$\pk c< \om^{\ve_{{a'}-1}+\pk 0}= \om^{\ve_{{a'}-1}+1}.$
We have $b<m< A_{a+1}0\leq A_{{a'}}0$ and the induction hypothesis yields
$\pk b< \om^{\ve_{{a'}-1}+\pk 0}= \om^{\ve_{{a'}-1}+1 } .$
It follows that $\ve_{a-1}+\pk b <  \om^{\ve_{{a'}-1}+1 } $, hence
$\pk m = \om^{\ve_{a-1}+\pk b}+\pk c< \om^{\ve_{{a'}-1}+1 } \leq \pk n.$
\smallskip

\item ($a={a'}$). Note that then $b\leq b'$. We consider several sub-cases.
\smallskip

\begin{enumerate}[label*= .\arabic*,wide, labelwidth=!, labelindent=0pt]
\item ($b<{b'}$).
The induction hypothesis yields $\pk b< \pk {b'}$.
Hence $\om^{\ve_{a-1}+\pk b}<\om^{\ve_{a-1}+\pk {b'}}$.
We have $c<A_ab$, and the subsidiary induction hypothesis yields
$\pk c< \om^{\ve_{a-1}+\pk b}<\om^{\ve_{a-1}+\pk {b'}}$.
Putting things together we see
$\pk m = \om^{\ve_{a-1}+\pk b} +\pk c< \om^{\ve_{a-1}+\pk {b'}}  \leq \pk n.$
\smallskip

\item ($b={b'}$). The inequality $m<n$ yields $c<{c'}$ and the induction hypothesis yields
$\pk c< \pk {c'}$. Hence 
$\pk m = \om^{\ve_{a-1}+\pk b} +\pk c<\om^{\ve_{a-1}+\pk b} +\pk {c'}=\pk n$. \qedhere
\end{enumerate}
\end{enumerate}
\endproof

As before, our ordinal assignment is invariant under base change.

\begin{lemma} 
$ m\bchb k{k+1} \bchb {k+1}\omega=\pk m$.
\end{lemma}

\proof
Write $A_xy$ for $A_x(k,y)$ and $B_xy$ for $A_x(k+1,y)$ and proceed by induction on $m$.
The assertion is clear for $m=0.$
Let $m\gnf A_a b+c$.
Then, $m\bchb k{k+1}  \gnf B_a( b\bchb k{k+1} ) +c\bchb k{k+1} $, and the induction hypothesis yields

\begin{eqnarray*}
\pskpe {m\bchb k{k+1} }&=&\pskpe {\big (B_a( b\bchb k{k+1} ) +c\bchb k{k+1}  \big )}\\
&=&\om^{\ve_{a-1}+\pskpe {b\bchb k{k+1} }} +\pskpe {c\bchb k{k+1} }\\
&=&\om^{\ve_{a-1}+\pk b} +\pk c =\pk m.
\end{eqnarray*}
\endproof

As we did for the first Goodstein process, we define
\[\bo _k(\ell) := \BG_k(\ell) \bchb {k+2}\omega. \]
The ordinals $\bo _k(\ell)$ are decreasing on $k$ provided they are non-zero, from which the termination of the process follows.

\begin{theorem} For all $\ell<\om$, there exists a $k<\om$ such that $\BG_k(\ell)=0.$
This is provable in ${\sf PA}+{\rm TI}(\ve_\om)$.
\end{theorem} 

\proof
By the previous lemmata,
\begin{eqnarray*}
\bo_{k+1}(\ell )&=& \BG_{k+1}(\ell)\bchb {k+3}\omega
=  (\BG _{k}(\ell)\bchkb-1)\bchb {k+3}\omega \\
&<&\BG_{k}(\ell)\bchkb\bchb {k+3}\omega 
=\BG_{k}(\ell)\bchb {k+2}\omega= \bo_k(\ell ).
\end{eqnarray*}
Since $( \bo_k(\ell ))_{k<\omega}$ cannot be an infinite decreasing sequence of ordinals, there must be some $k$ with $ \bo_k(\ell ) = 0$, yielding $\BG_k(\ell) = 0$.
\endproof

Now we are going to show that for every $\al<\varepsilon_\omega$,
${\sf PA}+{\rm TI}(\al)\not\vdash \forall \ell \exists k \ b_k(\ell )=0.$
In view of Proposition \ref{propMajorize}, the following technical lemma will be crucial for proving our main independence result for ${\sf ACA}'_0$.

\begin{lemma}\label{majB}
Given $k,m<\omega$ with $k\geq 2$,
\[\pskpe {(m\bchb k{k+1} -1)} \geq \pk m[k-1].\]
\end{lemma}

\proof
Write $A_xy $ for $A_x(k,y) $ and $B_xy $ for $A_a(k+1,y)$.
We prove the claim by induction on $m$.
Let $m\gnf A_ab +c.$ 
\smallskip

\begin{enumerate}[label*={\sc Case \arabic*},wide, labelwidth=!, labelindent=0pt]

\item ($c>0$).
Then the induction hypothesis and Lemma \ref{lemPkMon} yield

\begin{align*}
\pskpe {(m\bchkb-1)} &= \om^{\ve_{a-1}+\pskpe{b\bchkb}}+ \pskpe{(c\bchkb-1)}\\
&\geq  \om^{\ve_{a-1}+\pk b } + \pk c  [k-1]
=  \left (\om^{\ve_{a-1}+\pk b }  +\pk c \right )[k-1]\\
&=   \pk{\left ( A_ab  + c\right )}[k-1]
=  \pk m[k-1].
\end{align*}

\item ($c=0$).
We consider several sub-cases.

\smallskip

\begin{enumerate}[label*= .\arabic*,wide, labelwidth=!, labelindent=0pt]

\item\label{case31} ($a>0$ and $b>0$). 
The induction hypothesis yields

\begin{align*}
\pskpe{(m\bchkb-1)}&  = \pskpe{(B_a( b\bchkb)-1)}\\
&\geq   \pskpe{B_a(  b\bchkb -1)}\cdot k
=   \om^{\ve_{a-1}+\pskpe{(b\bchkb-1)}}\cdot k \\
&\geq   \om^{\ve_{a-1}+ \pk b[k-1]}\cdot k
\geq    \om^{\ve_{a-1}+ \pk b} [k-1]
=  \pk m [k-1],
\end{align*}
since $B_a( b\bchkb-1) $ is in $k+1$ normal form by Lemma \ref{NF} and Lemma \ref{pr}.
\smallskip

\item  ($a>0$ and $b=0$). Then, the induction hypothesis yields

\begin{eqnarray*}
\pskpe{(m\bchkb-1)}
&=&\pskpe{(B_a0-1)}
=\pskpe{(B_{a-1}^{k+1}0-1)}\\
&=&\pskpe{(B_{a-1}^{k-1} B^2_{a-1} 0-1)}\\
&\geq &\pskpe{ (B_{a-1}^{k-1}2 -1) } > \pskpe{ B_{a-1}^{k-1}1  } \\
&\geq&    \om_{k-1}(\varepsilon_{a-1}+1)
=    \ve_{a-1} [k-1]
\\
&=&   \pk{ A_a 0}[k-1]
=  \pk m [k-1],
\end{eqnarray*}
where in the last inequality we use that $B_{a-1}^{\ell}0$ is in $k+1$ normal form for $\ell\leq k$ by Lemma \ref{NF} and Lemma \ref{pr}, and an easy induction on $k$.
\smallskip

\item  ($a=0$ and $b>0$). Then the  induction hypothesis yields similarly as in \ref{case31}:

\begin{eqnarray*}
\pskpe{(m\bchkb-1)}&=&\pskpe{(B_0(b\bchkb) -1)} = \pskpe{((k+1)^{b\bchkb} -1)}\\
&\geq  &   \pskpe{(  (k+1)^{ b \bchkb -1}\cdot k)}\\
&=  &   \om^{\pskpe{(b \bchkb  -1)}}\cdot k 
\geq    \om^{ \pk b [k-1]}\cdot k
>  \pk m[k-1],
\end{eqnarray*}
since $(k+1)^{ b\bchkb -1}\cdot k $ is in extended $(k+1)$-normal form.
\smallskip

\item ($a=0$ and $b=0$).
Recall that by convention, $\ve_{-1} = 0$.
Then, $\pk m[k-1] = \pk{B_00}[k-1] = \omega^{\ve_{-1}+0} [k-1] = \omega^0[k-1] = 0  $, and the claim follows.\qedhere
\end{enumerate}

\end{enumerate}
\endproof

\begin{theorem}\label{ThmA}
${\sf ACA}'_0 \not \vdash \forall \ell \exists k \ \BG_k(\ell )=0.$
\end{theorem}

\proof 
The proof runs similarly to that of Theorem \ref{ThmTermA}.
This time, we define $\ell_n = A_{n+1}0$.
Observe that $\ell_n $ is in normal form and $\pk(\ell_n) = \omega^{\ve_n + 0} = \ve_n$.
Moreover, $  \ve_n = \ve_\om[n]$.
It follows in view of Lemma \ref{majB} that $\forall \ell \exists k \ \BG_k(\ell )=0$ implies that $F_{\ve_\om}$ is total, which is not provable in ${\sf PA}+{\rm TI}(\al)$ for any $\al<\ve_\om$.
\endproof

\section{Goodstein sequences for ${\sf ACA}_0^+$}

Finally, we consider what should intuitively be the strongest version of our Goodstein process: that where the base change is simultaneously applied to $a$, $b$, and $c$.
The resulting Goodstein principle will then be independent of ${\sf ACA}_0^+$.
The intuition should be that an expression $A_ab+c$ is written so that each of $a,b,c$ is hereditarily represented using the Ackermann function.

\begin{definition}\label{defBCHplus}
If $2\leq k<\ell<\omega$ and  $m\in \mathbb N$, define $m\bcc k\ell$ recursively by:
\begin{enumerate}
\item $0\bcc k\ell:=0$
\item $m\bcc k\ell:=A_{a\bcc k\ell}\left (\ell,b\bcc k\ell\right ) +c\bcc k\ell$ if $m\gnf A_a(k,b) +c$.
\end{enumerate}
We write $\bchkc$ instead of $\bcc k{k+1}$.
\end{definition}
As before, we observe that if $m\gnfm A_a(k,b)\cdot d+ e$, then
\[ m \bcc k\ell = A_{a\bcc k\ell}(\ell,b\bcc k\ell) \cdot d + e \bcc k\ell.\]
We can then define our final Goodstein process based on this new base change operator.

\begin{definition}
Let $\ell<\om$. 
Put $\CG_0(\ell):=\ell.$
Assume recursively that $\CG_k(\ell)$ is defined and $\CG_k(\ell)>0$.
Then, $\CG_{k+1}(\ell)=\CG_k(\ell)\bcc{k+2}{k+3}-1$. If $\CG_k(\ell)=0$, then $\CG_{k+1}(\ell):=0$.
\end{definition}

Termination and independence results can then be obtained following the same general strategy as before.

\begin{lemma}\label{lemUnMonAB}
If $ m < n $ and $k\geq 2$, then $ m \bchkc< n\bchkc$.
\end{lemma}

\proof
The assertion is clear if $m=0$. As before we write $A_xy$ for $A_x(k,y)$, $B_xy$ for $B_x(k+1,y)$, $m\gnf A_ab +c$
and $n \gnf A_{{a'}} {b'} + {c'}$.
\begin{enumerate}[label*={\sc Case \arabic*},wide, labelwidth=!, labelindent=0pt]
\item ($a < {a'}$).\label{case12}
By the induction hypothesis $a <{a'}$ yields $a\bchkc \leq {a'}\bchkc$.
Write $m\gnfm A_ab\cdot d + e$.
We have
\[m \bchkc= B_{a\bchkc}( b\bchkc)\cdot d +e\bchkc .\]
There are a few things to consider:
first, we have $A_ab\leq m < A_{a+1}0 = A_{a}^k0= A_a A_a^{k-1}0$. Then $b < A_a^{k-1}0$. By the induction hypothesis, $b\bchkc < B_{a\bchkc}^{k-1}0$.
We have that $d<A_{a+1}0$.
Moreover $A_ab$ is in $k$-normal form by lemma \ref{NF} and  $e<m\leq A_ab$ yields
$e\bchkc<B_{a\bchkc}(  b \bchkc) =  B^k_{a\bchkc} 0$ using our bound for $b \bchkc$.
This yields
\begin{align*}
m\bchkc&
= B_{a\bchkc }(  b\bchkc )\cdot d + e\bchkc\\
&<   B_{a\bchkc }^{k }0\cdot A_{a+1}0 + B_{a\bchkc }^{k }0\\
&\leq  ( B_{a\bchkc }^{k }0)^2 +  B_{a\bchkc }^{k }0\\
& \leq B^{k+1} _{a\bchkc } 0 =  B _{a\bchkc + 1 }0\\
&\leq B _{{a'}\bchkc   }({b'}\bchkc ) + {c'}\bchkc = n \bchkc,
\end{align*}
where we use Lemma \ref{lemmUseful}.\ref{itSquare} for the third inequality.

\smallskip \item ($a={a'}$ and $b<{b'}$).
Consider two sub-cases.
\smallskip

\begin{enumerate}[label*= .\arabic*,wide, labelwidth=!, labelindent=0pt]
\item ($A_a( b+1)<n$).
Since $m < A_a( b+1)$ and $b+1\leq {b'}$ then by the induction hypothesis
$m\bchkc<B_{a\bchkc}( (b+1)\bchkc) \leq B_{a\bchkc}( {b'}\bchkc) \leq n\bchkc$.

\smallskip \item ($A_a( b+1)=d$).
Write $m \gnfm A_ab\cdot d+e$.
Here we consider two cases, depending on the value of $a$.

\smallskip

\begin{enumerate}[label*= .\arabic*,wide, labelwidth=!, labelindent=0pt]
\item
($a=0$). In this case $m=A_0b\cdot d +e =k^b\cdot d+e<k^{b+1}=n$, where $d<k$,  $e<k^b$, and  $n$ has $k$-normal form $k^{b+1}$.
The induction hypothesis yields $b\bchkc<(b+1)\bchkc$ and $e\bchkc<(k+1)^{b\bchkc}$.
Then
\begin{align*}
m\bchkc & =(k+1)^{b\bchkc}\cdot d+e\bchkc<(k+1)^{b\bchkc } \cdot k + (k+1)^{b\bchkc} \\
& =(k+1)^{ b \bchkc + 1}
\leq (k+1)^{(b+1)\bchkc} =n\bchkc.
\end{align*}
\smallskip

\item($a>0$).
We have $e< A_ab$ and by the induction hypothesis,  $e \bchkc < B_{a\bchkc}( b\bchkc)$. Moreover $d < A_a(k,b+1)$.
Let $r= B_{a\bchkc-1}^kB_{a\bchkc}(b\bchkc)$, so that $B_{a\bchkc-1} r = B_{a\bchkc}(b\bchkc+1)$.
Note that $r\geq A^k_{a-1}A_ab = A_{a}( b+1)$.
Then,
\begin{align*}
m\bchkc
&=B_{a\bchkc}( b\bchkc  )\cdot d+e\bchkc\\
&<  B_{a\bchkc}( b\bchkc )   \cdot A_{a}( b+1)+B_{a\bchkc}( b\bchkc )\\
&\leq  r^2+r < B_{a\bchkc-1 } r = B_{a\bchkc}(b\bchkc+1)\\
&\leq B_{a\bchkc}(  {(b+1)}\bchkc) = n \bchkc
\end{align*}
where the third inequality uses Lemma \ref{lemmUseful}.\ref{itSquare}.
\end{enumerate}

\end{enumerate}


\smallskip \item ($a={a'}$ and ${b'}=b$). Since $m<n$ then it must be that $c<{c'}$, so
\begin{eqnarray*}
m\bchkc
&=&B_{a\bchkc}( b\bchkc  ) +c\bchkc\\
&<&B_{a\bchkc}( b\bchkc  ) +{c'}\bchkc
= n\bchkc .
\end{eqnarray*}
\end{enumerate}
\endproof

Thus, the base-change operation is monotone.
Next we see that it also preserves normal forms.

\begin{lemma}\label{prAB}
If $m=A_a(k,b) +c$ is in $k$-normal form, then
\[m\bchkc\gnfp A_{a\bchkc}(k+1,b\bchkc) +c\bchkc.\]
\end{lemma}

\proof
As usual write $A_xy$ for $A_x(k,y)$ and $B_xy$ for $A_x(k+1,y)$.
Let $m\gnfm A_ab\cdot d+e$.
We have that  $A_a0 \leq m<A_{a+1}0$ and $A_ab \leq m<A_a( b+1)$.
So, $A_ab < A_{a+1}0 = A_a^k 0 $. By Lemma \ref{NF}, $ A_a^{\ell} 0 $ is in $k$-normal form for $\ell<k$.
Hence by Lemma \ref{lemUnMonAB}, 
$b< A_a^{k-1}0 $ yields $b\bchkc< B_{a\bchkc}^{k-1}0 $.
We have that $d<m\leq A_a(b+1)$, and $A_a b$ is in $k$-normal form, Lemma \ref{lemUnMonAB} yields
$e \bchkc < B_{a\bchkc}( b\bchkc) < B_{a\bchkc}B_{a\bchkc}^{k-1}0 $.
Let $r=B_{a\bchkc}^{k}0$.
Then, 
\begin{align*}
m \bchkc  &= B_{a\bchkc} (b\bchkc)\cdot d + e\bchkc\\
&\leq  B_{a\bchkc}^{k}0\cdot A_{a\bchkc }^{k} 0  +  B_{a\bchkc}^{k}0\\
& \leq r^2+r \leq B_{a\bchkc +1}0.
\end{align*}
Now we check that $ c\bchkc < B_{a\bchkc}( b\bchkc + 1)$.

If $a=0$, then $c\gnfm A_0b\cdot d+e = k^b\cdot d+e$ with $d<k$ and $e<k^b$.
Then, $e\bchkc<(k+1)^{b\bchkc}$.
Thus
\begin{align*}
m\bchkc & =(k+1)^{b\bchkc}\cdot d+e\bchkc \\
&\leq (k+1)^{b\bchkc}\cdot k + (k+1)^{b\bchkc}= (k+1)^{b\bchkc+1}=B_0(b\bchkc+1),
\end{align*}
as needed.

In the remaining case, we have that $a>0$.
Let $r= B_{a\bchkc-1}^k B_{a\bchkc} (b\bchkc)$.
Then,
\begin{align*}
m\bchkc &
=B_{a\bchkc}( b\bchkc  )\cdot d+e\bchkc\\
&<B_{a\bchkc}( b\bchkc  )\cdot A_{a\bchkc}( b+1)  + B_{a\bchkc}(  b\bchkc  )\\
&\leq r^2+r
 \leq B_{a\bchkc}( b\bchkc+1  ),\\
\end{align*}
where the last inequality uses Lemma \ref{lemmUseful}.\ref{itSquare}.
So $B_{a\bchkc}( b\bchkc) +c\bchkc$ is in $(k+1)$-normal form.
\endproof

Finally, we provide an ordinal mapping to show that the Goodstein process terminates.
Recall that $\ve_{-1} = 0$ by convention, and $-1+\alpha$ is the unique $\eta$ such that $1+\eta=\alpha$, provided $\alpha>0$.
For the sake of legibility, we write $\epfun\alpha$ instead of $\ve_{-1+\alpha}$.

\begin{definition}\label{defOrdAB}
Given $k\geq 2$, define a function $\ck\cdot \colon \mathbb N \to \varphi_2(0)$ given by:
\begin{enumerate}
\item $\ck0:=0$.
\item $\ck m:=\om^{\epfun{\ck a}+\ck b}+\pk c$ if $m\gnf A_a(k,b) +c$.
\end{enumerate}
\end{definition}

As was the case for the previous mappings, the mappings in Definition \ref{defOrdAB} are strictly increasing and invariant under base change, as we show in the following lemmas.

\begin{lemma}\label{lemPkMonAB}
If $m<n $, then $\ck m<\ck n$.
\end{lemma}

\proof
Proof by induction on $d$ with subsidiary induction on $c$.
The assertion is clear if $c=0$.
Let $m\gnf A_ab  +c$
and $n\gnf A_{{a'}} {b'} +{c'}$.

\begin{enumerate}[label*={\sc Case \arabic*},wide, labelwidth=!, labelindent=0pt]
\item ($a<{a'}$). By the induction hypothesis $\ck a<\ck {a'}$ , so ${\epfun{\ck a}}<{\epfun{\ck {a'}}}.$
We have $c<m<A_{a+1}0\leq A_{{a'}}0\leq n$, so the induction hypothesis yields
\[\ck c< \om^{\epfun{\ck {a'}}+\ck 0}=\epfun{\ck {a'}}.\]
We have $b<m< A_{a+1}0\leq A_{{a'}}0$ and the induction hypothesis also yields
\[\ck b< \om^{\epfun{\ck {a'}}+\ck 0}=\epfun{\ck {a'}}.\]
It follows that $\epfun{\ck a}+\ck b <  \epfun{\ck {a'}}$, hence
\[\ck m = \om^{\epfun{\ck a}+\ck b} +\ck c<\epfun{\ck {a'}}\leq \ck n.\]
\smallskip


\item ($a={a'}$). We consider several sub-cases.
\smallskip

\begin{enumerate}[label*= .\arabic*,wide, labelwidth=!, labelindent=0pt]
\item ($b<{b'}$). The induction hypothesis yields $\ck b< \ck {b'}$.
Hence $\om^{\epfun{\ck a}+\ck b}<\om^{\epfun{\ck a}+\ck {b'}}$.
We have $c<A_a(k,b)$, and the subsidiary induction hypothesis yields
\[\ck c < \om^{\epfun{\ck a}+\ck b}<\om^{\epfun{\ck a}+\ck {b'}},\]
so
\[\ck m = \om^{\epfun{\ck a}+\ck b} +\ck c< \om^{\epfun{\ck a}+\ck {b'}}  \leq \ck n.
\]
\smallskip


\item ($b={b'}$). The inequality $m<n$ yields $c<{c'}$ and the induction hypothesis yields
$\ck c< \ck {c'}$. Hence 
\begin{align*}
\ck m & = \om^{\epfun{\ck a}+\ck b}+\ck c\\
&<\om^{\epfun{\ck a}+\ck b} +\ck {c'}=\ck n.\qedhere
\end{align*} 
\end{enumerate}
\end{enumerate}
\endproof

Next we prove that the ordinal assignment is invariant under base change.

\begin{lemma} 
For all $m<\omega$, $ m\bcc k{k+1}\bcc {k+1}\omega=\ck m$.
\end{lemma}

\proof
Proof by induction on $m$.
The assertion is clear for $m=0.$
Assume $c\gnf A_ab+c$.
Then, $c\bchkc \gnf B_{a\bchkc}( b\bchkc) +c\bchkc$, and the induction hypothesis yields

\begin{eqnarray*}
\cskpe {m\bchkc}&=&\cskpe {(B_{a\bchkc}( b\bchkc) +c\bchkc)}\\
&=&\om^{\epfun{\cskpe {a\bchkc}}+\cskpe {b\bchkc}} +\cskpe {c\bchkc}\\
&=&\om^{\epfun{\ck a}+\ck b} +\ck c =\ck m.\qedhere
\end{eqnarray*}
\endproof

With this, we define for $m,k\in\mathbb N$ 
\[\co_k(m) = {\CG_k (m)}\bcc {k+2}\omega.\]
As in the proof of e.g.~Theorem \ref{ThmTermA}, the sequence $(\co_k(m))_{k<\omega}$ is decreasing on $k$ as long as $\CG_k (m) > 0$.
Since $\co_k(m)\leq \varphi_20$, we obtain the following.

\begin{theorem}
For all $\ell <\om$, there exists a $k<\om$ such that $c_k(\ell )=0.$ This is provable in ${\sf PA}+{\rm TI}(\varphi_20)$.
\end{theorem} 

Finally, we show that for every $\al<\varphi_20$,
${\sf PA}+{\rm TI}(\al)\not\vdash \forall \ell \exists k \ c_k(\ell )=0.$
For this, we need the following analogue of Lemma \ref{majA}.

\begin{lemma}\label{majAB}
Given $k,m<\omega$ with $k\geq 2$, $\cskpe {(m\bchkc-1)} \geq  \ck m [k-1]$.
\end{lemma}

\proof
We prove the claim by induction on $m$.
Write $A_xy$ for $A_x(k,y)$ and $B_xy$ for $A_x(k+1,y)$.
Let $m\gnf A_ab+c \gnfm A_ab\cdot d+e$. 
\smallskip

\begin{enumerate}[label*={\sc Case \arabic*},wide, labelwidth=!, labelindent=0pt]

\item ($c>0$).
Then the induction hypothesis and Lemma \ref{lemPkMonAB} yields
\begin{align*}
\cskpe{(m\bchkc-1)} &= \cskpe {( B_{a\bchkc}( b\bchkc) +c\bchkc-1)}\\
&= \om^{\epfun{\cskpe{a\bchkc}}+\cskpe{b\bchkc}}  + \cskpe{(c\bchkc-1)}\\
&\geq  \om^{\epfun{\ck a}+\ck b}  + \ck c[k-1]\\
&=  (\om^{\epfun{\ck a}+\ck b}  +\ck c)[k-1]\\
&=   \ck{( A_ab + c)}[k-1]
= \ck c[k-1].
\end{align*}

\item ($c=0$).
We consider several sub-cases.
\smallskip

\begin{enumerate}[label*= .\arabic*,wide, labelwidth=!, labelindent=0pt]

\item\label{case31} ($a>0$ and $b>0$). 
We have $m=A_ab$. The induction hypothesis yields
\begin{align*}
\cskpe{(m\bchkc-1)}&  = \cskpe{(B_{a\bchkc}( b\bchkc)-1)}\\
&\geq \cskpe{(B_{a\bchkc}(  b\bchkc -1) \cdot k)}\\
&= \left ( \om^{\epfun{\cskpe a\bchkc}+\cskpe{(b\bchkc-1)}}\right ) \cdot k\\
&\geq \left (  \om^{\epfun{\ck a}+ \ck b[k-1]} \right) \cdot k \\
&>  \left ( \om^{\epfun{\ck a}+ \ck b }\right )[k-1]
=  \ck m[k-1].
\end{align*}
\smallskip

\item  ($a>0$ and $b=0$). 
In this case $m=A_a0$. Then, the induction hypothesis yields
\begin{align*}
\cskpe{(m\bchkc-1)}& =\cskpe{(B_{a\bchk}0-1)}\\
&= \cskpe{(B_{a\bchkc -1}^{k+1}0-1)} \\
&\geq  \cskpe{(B_{a\bchkc -1}^{k-1}1)}\\
&\geq  \om_{k-1}\left (\epfun{\cskpe{(a\bchkc -1)}}+1 \right ) \\
& \geq   \om_{k-1}\left (\epfun{\ck a [k-1]}+1 \right )\\
&\geq     \epfun{\ck a} [k-1]
=  \ck{ A_a0 }[k-1],
\end{align*}
Since $B_{a\bchkc -1}^{\ell}1$ is in $(k+1)$-normal form for $\ell<k+1$.
\smallskip

\item  ($a=0$ and $b>0$). Then the  induction hypothesis yields
\begin{eqnarray*}
\cskpe{(m\bchkc-1)}&=&\cskpe{(B_0(b\bchkc)-1)}\\
&\geq&\cskpe{(B_0(b\bchkc-1)\cdot k)}\\
& = &     \cskpe{ (k+1)^{(b\bchkc-1)}}\cdot k\\
&=  &   \om^{\cskpe{(b\bchkc-1)}} \cdot k\\ 
&\geq &    \om^{ \ck b [k-1]}\cdot k
\geq  \ck m[k-1],
\end{eqnarray*}
since $(k+1)^{ b\bchkc-1 }\cdot k $ is in $(k+1)$-normal form.
\smallskip

\item ($a=0$ and $b=0$).
In this case $\ck m[k-1] = \ck {A_00}[k-1] = \omega^0[k-1]=0$, so the lemma follows.

\end{enumerate}

\end{enumerate}
\endproof

\begin{theorem}
${\sf ACA}_0^+\not\vdash \forall \ell \exists k \CG_k(\ell )=0.$
\end{theorem}
\proof 
The proof once again runs similarly to that of Theorem \ref{ThmTermA} using in this case Lemma \ref{majAB}.
In this case, we define $\ell_n$ recursively by $\ell_0 = 1$ and $\ell_{n+1} = A_{\ell_n}(2,0)$, and observe that $\ck{\ell_{n+1}}= \varphi^n_1 1 =  {\varphi_2 0}[n] $.
\endproof

\section{Concluding remarks}

This work follows on Arai et al.~\cite{FSGoodstein}, where it is shown that a more elaborate normal form based on the Ackermann function leads to an independence result of strength $\Gamma_0$.
The simpler normal forms considered here, which are arguably more natural, lead to substantially weaker Goodstein principles.
However, it is instructive to observe that changes in the notion of base-change can lead to somewhat unpredictable modifications in proof-theoretic strength, naturally leading to independence results for ${\sf ACA}_0$, ${\sf ACA}'_0$, and ${\sf ACA}^+_0$, three prominent theories of reverse mathematics.

A natural question is whether the gap between $\ve_0$ and $\Gamma_0$ can be bridged using further variants of the Ackermannian Goodstein  process.
One idea would be to perform the sandwiching of \cite{FSGoodstein} only up to a fixed number $n$ of steps.
One may conjecture that such a process would lead to independence results for ${\sf TI}(\Gamma_0[n])$, where the fundamental sequence for $\Gamma_0$ is defined by iteration on the parameter $\alpha$ in $\varphi_\alpha\beta$.
However, the results presented here should serve as evidence that the proof-theoretic strength of Goodstein principles is not easy to predict, so such conjectures should be taken with a grain of salt!

\subsection*{Acknowledgements}

This work was partially supported by the FWO-FWF Lead Agency grant G030620N (FWO)\slash I4513N (FWF).

\ignore{

}

 \bibliographystyle{plain}
 \bibliography{biblio}

\end{document}